\documentclass[a4paper,11pt]{article}
\usepackage{array}

\usepackage{amsmath,amssymb} 
\usepackage{amsthm}
\usepackage{float} 
\input xy
\xyoption{all}
\usepackage[dvips]{color}

\newcommand{\bfalpha}{\mbox{\boldmath$\alpha$}}

\newcommand{\comment}[1]{}
\newcommand{\PP}{{\mathbb{P}}}

\newcommand{\CC}{{\mathbb{C}}}
\newcommand{\ZZ}{{\mathbb{Z}}}
\newcommand{\eps}{\varepsilon}
\newcommand{\arl}{\ar@{..>}[l]}
\newcommand{\ard}{\ar@{..>}[d]}

\usepackage[dvips]{graphicx}
\usepackage{rotating}

\newtheorem{theorem}{Theorem}

\newtheorem{lemma}[theorem]{Lemma}
\newtheorem{corollary}[theorem]{Corollary}

\theoremstyle{definition}

\newtheorem{remark}[theorem]{Remark}
\newtheorem{example}[theorem]{Example}

\begin{document}
\title{
Sums of squares of binomial coefficients,
with applications to Picard-Fuchs equations}

\author{H. A. Verrill}

\maketitle

\begin{abstract}
For an arbitrary fixed positive integer $N$, we give a 
recurrence relation (\ref{recrel2})  for the 
sequence
$a_n^N:=\sum\left(\frac{n!}{p_1!p_2!\cdots p_N!}\right)^2$,
where the sum runs over sets of non-negative integers $p_1,\dots,p_N$
summing to $n$.
For the more general sequence
$a_n^{\alpha_1,\dots,\alpha_N}
:=\sum
\left(\prod_{i=1}^N
\alpha_i^{p_i}\right)
\left(\frac{n!}{p_1!p_2!\cdots p_N!}\right)^2,$ for any fixed $\alpha_1, \dots
 \alpha_N \in \CC$,
we give an algorithm for determining a recurrence relation.  
We show how to apply the latter procedure to obtain general 
recurrence relations 
((\ref{recrel_generalN=2case}) and (\ref{recrelN=3abc}))
for the cases $N=2$ and $3$.  
As a geometric application,
in Section~\ref{sec:ellipticPF}
 we give the Picard-Fuchs equation 
(\ref{ellipticPFeqn})
for the families of elliptic curves corresponding to the $N=3$ case.

The main idea is to introduce auxiliary sequences 
(\ref{eqn_aux_terms}) or (\ref{eqn_aux_terms_general}), 
generalising and systematising the method used in 
\cite{verrillKyoto} for the case $N=3, \alpha_1,\alpha_2,\alpha_3=1$.
This approach seems to be similar to the method used by Cusick 
\cite{C}
in finding recurrences for sums of powers of binomials.
We briefly discuss the application to this case in 
Section~\ref{sec:furtherapplications}.
\end{abstract}

\section{Introduction}

In this paper we study sequences $\displaystyle{\{a_n^{\bfalpha}\raisebox{0cm}[0cm][0.11cm]{\}}_n}$ where
\begin{equation}
\label{eqn:firstdefinition}
a_n^{\{\alpha_1,\dots,\alpha_N\}}
:=
\sum_{p_1+\cdots+p_N=n}
\left(\prod_{i=1}^N
\alpha_i^{p_i}\right)
\left(\frac{n!}{p_1!p_2!\cdots p_N!}\right)^2,
\end{equation}
$N\ge 2$,
${\bfalpha}=
\{\alpha_1,\dots,\alpha_N\}\in\CC^N$ is fixed,
and in the sum
$p_i$ are always non-negative integers.
We will use the notation 
$\binom{n}{p_1,\dots,p_N}=\left(\frac{n!}{p_1!p_2!\cdots p_N!}\right)$,
write 
$a_n^{\alpha_1,\dots,\alpha_N}=a_n^{\{\alpha_1,\dots,\alpha_N\}}$,
and use the notation $a_n^N:=a_n^{1,\dots,1}$, where the superscript
is a list of $N$ ones.  This should cause no confusion, since we
do not allow $N=1$ in definition 
(\ref{eqn:firstdefinition}), as that case is too easy.

The simplest non-trivial example is the sequence $a_n^{1,1}$, for which
we have the following sequence
(\#A000984 in Sloan's tables \cite{sloan}).
\begin{equation}
a_n^2:=a_n^{1,1}:=\sum_{p=0}^n\binom{n}{p}^2=\binom{2n}{n}.
\label{ex:simplestcase}
\end{equation}
This is a well known result, a special case
of $\sum_{p+q=C}\binom{A}{p}\binom{B}{q}=\binom{A+B}{C}$,
the number of ways of taking $C$ elements from a set which is the
union of two sets of sizes $A$ and $B$, by summing over the number of
elements taken from each set.

%(\ref{recrel_generalN=2case})
Another example, similar in appearance,
 is the case of $a^{1,-1}_n$, where
\begin{equation}
\label{motivating_1}
a^{1,-1}_n=\sum_{p=0}^n(-1)^{n-p}\binom{2n}{p}^2=\binom{2n}{n}.
\end{equation}
Examples (\ref{ex:simplestcase}) and (\ref{motivating_1})
may be found in many texts on combinatorics, and also in the
survey \cite{W}.
Along these lines we also have the less obvious result
\begin{equation}
\label{motivating_2}
a_n^{1,\omega,\omega^2}:=\sum_{p+q+r=3m}\omega^{p-q}\binom{3m}{p,q,r}^2
=\binom{4m}{2m,m,m},
\end{equation}
where $\omega$ is a primative third root of unity.
Proofs of (\ref{motivating_1}) and (\ref{motivating_2}) 
(given in Examples~\ref{ex:11_1-1} and \ref{ex:omega3} respectively)
will follow
from the 
general
recurrence relations we find, for example, Theorem~\ref{thm_ab_case}, which
says that in the $N=2$ case we have
\begin{equation}
na_n^{a,b} - (2n-1)(a+b)a_{n-1}^{a,b}
 + (n-1)(a-b)^2a_{n-2}^{a,b}=0.
\end{equation}
Theorem~\ref{thm_abc_case} gives a similar formula for the case $N=3$.
Corollary~\ref{cor:mainresult} states that we can always obtain 
recurrence relations, with at most $2^{N-1}+1$ terms, 
and Section~\ref{weightedsumsmethod} gives a concrete method of finding
these, which does not make use of computing any terms, and thus is
faster than methods which find recurrences given only the existance
of a recurrence of known degree, and the terms of the sequence.
Theorem~\ref{thm:recrel_for_N=1case} gives an explicit formula
for the recurrence for all $N$, in the case that all $\alpha_i$ are $1$,
and a few examples obtained by plugging in given values of $N$ are
given in Table~\ref{table:recrela=1}.

\subsection*{Geometrical motivation}

The motivation for studying the expressions given in the title of this note
comes from the study of the family of Calabi-Yau 
varieties, the general member of which is
given by the resolution of a hypersurface in $\PP^{N-1}$ defined by the
following equation---which actually should be multiplied through on both sides
by $X_1\dots X_N$ to obtain a homogeneous equation of degree $N-1$.
\begin{eqnarray}
\mathcal X^{\alpha_1,\dots,\alpha_N}_t:\>\>\>
(X_1+\cdots + X_N)\left(\frac{\alpha_1}{X_1}+\cdots + \frac{\alpha_N}{X_N}\right)t=1.
\label{eqn_for_A_{n-1}}
\end{eqnarray}
Here $t$ is the parameter of the family, and we take
$\alpha_i$ to be nonzero fixed elements of $\CC$.
Aspects of the $N=3$ case are considered in \cite{verrillKyoto} and \cite{verrillmodularity}, 
and the $N=4$ case is studied in \cite{HV}.  Some work on the general case can be
found in \cite{Ludwig}.

In order to study a family of varieties, one often works with the periods of the
members.  Using the methods of \cite{PetersStienstra}, in the $N=3,\alpha_i=1$
case, it was shown in \cite{verrillKyoto} that one of the periods of this
family has the form {\small{\it constant}} $\times\sum_{n\ge 0} a_n^{1,1,1} t^n$.  Exactly the same considerations
show that in general, up to a constant factor,
$\sum_{n\ge 0} a_n^{\alpha_1,\dots,\alpha_N} t^n$ is a period for 
$\mathcal X^{\alpha_1,\dots,\alpha_N}_t$. (See abstract or below for definitions of the
$a_n$.)

The Picard-Fuchs equation is a differential equation satisfied by the periods.
One method to find this equation (used for example in \cite{PetersStienstra})
is to find a recurrence relation for the coefficients $a_n$, which is what
the rest of this note is devoted to.

The rest of this paper is more combinatorial than geometric. 
We return to the geometry in Section~\ref{sec:ellipticPF}, where we
give the Picard-Fuchs equation for the families of elliptic curves, corresponding to the
$N=3$ case (\ref{ellipticPFeqn}).
However, we do not discuss the geometry or modularity of the varieties $\mathcal X^{\alpha_1,\dots,\alpha_N}_t$
for $N>3$ here, which is currently work in progress.
We do give differential equations (\ref{diffeqn1})
closely related to the Picard-Fuchs equations
in the case where all $\alpha_i=1$;  Table~\ref{tableofPFensa=1} gives some examples.

\section{The case $\alpha_1=\cdots=\alpha_N=1$}

\begin{table}
\begin{eqnarray*}
0 &=&na_n^2 - 2(2n -1)a_{n-1}^2\\
0 &=&n^2a_n^3 - (10n^2 - 10n + 3)a^3_{n-1} + 9(n-1)^2a^3_{n-2}\\
0 &=&n^3a_n^4 -2(2n-1)(5n^2 -5n +2)a_{n-1}^4 + 64(n-1)^3a_{n-2}^4\\
0 &=&n^4a_n^5 - (35n^4 - 70n^3 + 63n^2 - 28n + 5)a^5_{n-1}\\
&&+ (n-1)^2(259(n-1)^2 + 26)a_{n-2}^5 - (3\cdot 5)^2(n-1)^2(n-2)^2a_{n-3}^5\\
0 &=& n^5  a_n^6
      -2(2n-1)(14n^4 - 28n^3 + 28n^2 - 14n + 3)  a_{n-1}^6\\
&&       +4(n-1)^3(196(n-1)^2 + 59) a_{n-2}^6\\
&&      -(2\cdot4\cdot6)^2(n-1)^2(n-2)^2(n-
\mbox{$\frac{1}{2}$})^2 a_{n-3}^6\\
0 &=& n^6  a_n^7
      -(84(n(n-1))^3+126(n(n-1))^2 + 54n(n-1)+7)  a_{n-1}^7\\
&&       +3(n-1)^2( 658(n-1)^4 +  396(n-1)^2 + 17) a_{n-2}^7\\
&&       -2(n-1)^2(n-2)^2(6458n^2 - 19374n + 15505) a_{n-3}^7\\
&&        +(3\cdot5\cdot7)^2(n-1)^2(n-2)^2(n-3)^2a_{n-4}^7\\
0 &=& n^7  a_n^8\\
&&       -2(2n - 1)(30(n(n-1))^3+54(n(n-1))^2 + 27n(n-1)+4) a_{n-1}^8\\
&&       +12(n-1)^3(364(n-1)^4+365(n-1)^2+47) a_{n-2}^8\\
&&       -2^7(n-1)^2(n-2)^2(2n-3)(205n^2 - 615n + 554) a_{n-3}^8\\
&&        +(2\cdot4\cdot6\cdot8)^2(n-1)^2(n-2)^3(n-3)^2a_{n-4}^8\\
0 &=& n^8a_n^9\\
&&-(165n^8 - 660n^7 + 1386n^6 - 1848n^5 + 1650n^4 - 990n^3 + 385n^2 - 88n + 9)a_{n-1}^9\\
&&+3(n - 1)^2(2926n^6 - 17556n^5 + 48290n^4 - 76120n^3 + 71423n^2 - 37422n + 8487)a_{n-2}^9\\
&&-(n - 2)^2(n - 1)^2(172810n^4 - 1036860n^3 + 2489234n^2 - 2801832n + 1237167)a_{n-3}^9\\
&&+9(n - 3)^2(n - 2)^2(n - 1)^2(117469n^2 - 469876n + 493542)a_{n-4}^9\\
&&-(3\cdot5\cdot7\cdot9)^2
(n - 4)^2(n - 3)^2(n - 2)^2(n - 1)^2a_{n-5}^9\\
0&=&n^9a_{n}^{10}\\
&&-2(2n - 1)(55n^8 - 220n^7 + 484n^6 - 682n^5 + 649n^4 - 418n^3 + 176n^2 - 44n + 5)a_{n-1}^{10}\\
&&4(n - 1)^3(4092n^6 - 24552n^5 + 69993n^4 - 116292n^3 + 116754n^2 - 66396n + 16675)a_{n-2}^{10}\\
&&-8(n - 2)^2(n - 1)^2(2n - 3)(30580n^4 - 183480n^3 + 458909n^2 - 551067n + 267900)a_{n-3}^{10}\\
&&256(n - 3)^2(n - 2)^3(n - 1)^2(21076n^2 - 84304n + 97035)a_{n-4}^{10}\\
&&-(2\cdot4\cdot6\cdot8\cdot10)^2
(n - 4)^2(n - 3)^2(n - 2)^2(n - 1)^2(n - \mbox{$\frac{5}{2}$})a_{n-5}^{10}
%\\
%
%0&=&n^{10}a_{n}^{11}\\
%&&-(286n^10 - 1430n^9 + 3861n^8 - 6864n^7 + 8580n^6 - 7722n^5 + 5005n^4 - 2288n^3 + 702n^2 - 130n + 11)n^{10}a_{n-1}^{11}\\
%&&(n - 1)^2(28743n^8 - 229944n^7 + 885456n^6 - 2093520n^5 + 3273699n^4 - 3430284n^3 + 2333734n^2 - 936884n + 169125)n^{10}a_{n-2}^{11}\\
%&&-2(n - 2)^2(n - 1)^2(617474n^6 - 5557266n^5 + 22234927n^4 - 50050572n^3 + 66294111n^2 - 48677850n + 15401892)n^{10}a_{n-3}^{11}\\
%&&(n - 3)^2(n - 2)^2(n - 1)^2(21967231n^4 - 175737848n^3 + 553739264n^2 - 809054272n + 460053891)n^{10}a_{n-4}^{11}\\
%&&-9(n - 4)^2(n - 3)^2(n - 2)^2(n - 1)^2(14312974n^2 - 71564870n + 93063861)
%n^{10}a_{n-5}^{11}\\
%&&-(3\cdot5\cdot7\cdot9\cdot11)^2
%(n - 5)^2(n - 4)^2(n - 3)^2(n - 2)^2(n - 1)^2n^{10}a_{n-6}^{11}\\
\end{eqnarray*}
\caption{Recurrence relations for $a_n^N$, 
from Equation (\ref{recrel2}) for
$2\le N\le 10$}
\label{table:recrela=1}
\end{table}

We first consider the case $\alpha_1=\cdots\alpha_N=1$, since this is simpler, and in this
case we can obtain an explicit formula for the recurrence relation.
In order to find a recurrence relation for
$$a_n^N:=\sum_{\sum p_i=n}\left(\frac{n!}{p_1!p_2!\cdots p_N!}\right)^2,$$
we introduce auxiliary sequences of numbers $a_n^{N,j}$, 
for $0\le j\le N$,
defined by
\begin{equation}
a_n^{N,j}:=
\sum_{\sum p_i=n}\left(\frac{n!}{p_1!p_2!\cdots p_N!}\right)^2 
p_1\cdots p_j.
\label{eqn_aux_terms}
\end{equation}
We have $a_n^N=a_n^{N,0}$, and for convenience, we set
$a_n^{N,N+1}=a_n^{N,-1}=0$.

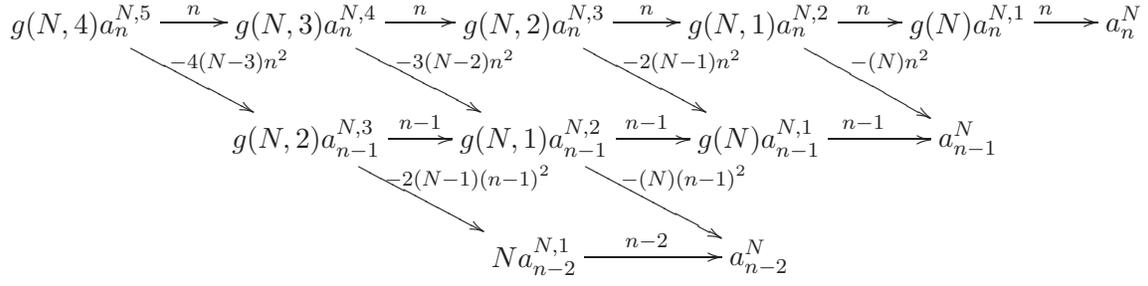
\begin{figure}
$$\hspace{-1cm}
\xymatrix{
g(N,4)a_n^{N,5}\ar[r]^n
\ar[dr]^{-4(N-3)n^2}
 & 
g(N,3)a_n^{N,4}\ar[r]^n
\ar[dr]^{-3(N-2)n^2}
 & 
g(N,2)a_n^{N,3}\ar[r]^n
\ar[dr]^{-2(N-1)n^2}
 & 
g(N,1)a_n^{N,2}\ar[r]^n
\ar[dr]^{-(N)n^2}
 & 
g(N)a_n^{N,1}\ar[r]^n
 & 
    a_n^{N}
\\
&
g(N,2)a_{n-1}^{N,3}\ar[r]^{n-1}
\ar[dr]^{-2(N-1){(n-1)}^2}
 & 
g(N,1)a_{n-1}^{N,2}\ar[r]^{n-1}
\ar[dr]^{-(N){(n-1)}^2}
 & 
g(N)a_{n-1}^{N,1}\ar[r]^{n-1}
& 
    a_{n-1}^{N}
\\
&
&
Na_{n-2}^{N,1}\ar[r]^{n-2}
 & 
a_{n-2}^{N}
}
$$
\caption{Diagram of relations (\ref{eqn_aux_reln_reversed}) between $a_n^{N,j}$ (\ref{eqn_aux_terms})}
\label{fig_arrows1}
\end{figure}

It is easy to verify that for $0\le j\le N$,
\begin{equation}
na_n^{N,j}=(N-j)a_n^{j+1} + jn^2a_{n-1}^{N,j-1}.
\label{eqn_aux_reln}
\end{equation}
Rearranging this equation, we have
\begin{equation}
(N-j)a_n^{N,j+1} =
na_n^{N,j}- jn^2a_{n-1}^{N,j-1}.
\label{eqn_aux_reln_reversed}
\end{equation}
Now, starting with $j=N$, and applying this relation repeatedly,
we will obtain a recurrence relation for the $a_n$.  We can visualise
this process in Figure~\ref{fig_arrows1},
drawn up to $a_n^{N,5}$,
where we use the notation
$$g(N,j)=(N-1)(N-2)\cdots(N-j).$$
In Figure~\ref{fig_arrows1}, the configuration 
\begin{equation}
\raisebox{4.2ex}{\mbox{
\xymatrix{
X \ar[r]^a \ar[rd]^b & A\\
& B
}}}
\label{eqn_meaningofdiagram}
\end{equation}
means $X = aA + bB$.
Thus to find an expression for
$g(N,4)a_n^{N,5}$ in terms of the $a_n^N$, we just have to follow
all paths from 
$g(N,4)a_n^{N,5}$ to the $a_n^N, a_{n-1}^N, a_{n-2}^N$,
multiplying the coefficients on the arrows in any path, and summing
the products for paths ending at the same point.
In this example, we obtain
\begin{eqnarray*}
g(N,4) &=&n^5a_n^N \\
&&- \Big(Nn^5 + 2(N-1)n^4(n-1) \\
&& \hspace{2cm} + 3(N-2)n^3(n-1)^2 + 4(N-3)n^2(n-1)^3\Big)
a_{n-1}^N 
\\&&+
\Big(3(N-2)Nn(n-1)^2 
+4N(N-3)(n-1)^3\\
&&\hspace{3.7cm} + 8(N-3)(N-1)(n-1)^2(n-2)\Big)a_{n-2}^N.
\end{eqnarray*}
This is true for all $N\ge 1$.  However, for $N=4$,
we have
$g(4,4)=0$, and so, substituting $N=4$,
and dividing by $n^2$, we obtain the relation:
$$n^3a_n^4 
-2(n-1)(5n^2 - 5n + 2)a_{n-1}^4
+64(n-3)^3 a_{n-2}^4=0.$$

The first few cases of $g(N,j)a_n^{N,j+1}$ are as follows.
\begin{eqnarray*}
g(N,0)a_n^{N,1} &=&na_n^{N} \\
g(N,1)a_n^{N,2} &=&n^2a_n^{N} - Nn^2a^N_{n-1} \\
g(N,2)a_n^{N,3} &=&n^3a_n^{N} -
\big[(3N - 2)n - (2N - 2)\big]n^2a^N_{n-1} \\
g(N,3)a_n^{N,4} &=&n^4a_n^{N} - 
[(6N - 8)n^2 + (-8N + 14)n + 3(N - 2)]n^2a^N_{n-1} \\
&& + 3N(N-2)(n-1)^2n^2a^N_{n-2}
\end{eqnarray*}

In general, the formula we obtain is:
\begin{equation}
g(N,j)a_n^{N,j+1}
=
\sum_{k\ge 0}
\left[
n^{j+1}
\!\!\!\!\!\!\!\!\!\!\!\!\!\!\!
\sum_{\scriptsize{\begin{array}{c}
1\le i\le k\\
\alpha_i+\beta_i=N+1,\alpha_i\in\mathbb N\\
1< \alpha_{i+1}+1< \alpha_i\le j
\end{array}
}}
\!\!\!\!\!\!\!\!\!\!\!\!\!\!\!
\prod_{i=1}^k
-\alpha_i\beta_i\left(\frac{n-i}{n-i+1}\right)^{\alpha_i-1}
\right]a_{n-k}^N,
\label{eqn_for_recrel}
\end{equation}
where an empty product is taken to be $1$.
Note that this is a finite sum, since for $k> \lfloor{N/2}\rfloor$,
there are no possible sequences of $\alpha_1,\dots,\alpha_k$ with 
$1\le \alpha_k< \alpha_1\le  j$ and $\alpha_{i+1}\le \alpha_i-2$.
Note also that since $\alpha_i>\alpha_{i+1}$ this is in fact a polynomial
in $n$.

Since $g(N,N)=0$, substituting $j=N$ in (\ref{eqn_for_recrel}), we obtain the
following result.
\begin{theorem}
\label{thm:recrel_for_N=1case}
For any positive integer $N$, the sequence $a_n^N
=\sum\left(\frac{n!}{p_1!p_2!\cdots p_N!}\right)^2$ (where
the sum is over non-negative integers $p_1,\dots,p_N$ summing to $n$)
satisfies the recurrence relation
\begin{equation}
\sum_{k\ge 0}
\left[
n^{N+1}
\!\!\!\!\!\!\!\!\!\!\!\!\!\!\!
\sum_{\scriptsize{\begin{array}{c}
1\le i\le k\\
\alpha_i+\beta_i=N+1,\alpha_i\in\mathbb N\\
1< \alpha_{i+1}+1< \alpha_i\le N
\end{array}
}}
\!\!\!\!\!\!\!\!\!\!\!\!\!\!\!
\prod_{i=1}^k
-\alpha_i\beta_i\left(\frac{n-i}{n-i+1}\right)^{\alpha_i-1}
\right]a_{n-k}^N
=0.
\label{recrel2}
\end{equation}
\end{theorem}
The recurrences in the first few cases, (after dividing by $n^2$), are given in Table~\ref{table:recrela=1}.

\begin{remark}
The coefficients of $a^N_{n-k}$ in (\ref{recrel2})
 may be written in other forms, for
example, the coefficient of $a^N_{n-1}$ can also be written as:
\begin{equation}
N\left(n^{N+2}-(n-1)^{N+2}\right)
-(N+2)n(n-1)\left(n^{N}-(n-1)^{N}\right).
\end{equation}
\end{remark}

\begin{remark}
Formula (\ref{recrel2}) is a closed formula.
However, alternatively, using (\ref{eqn_aux_reln_reversed}), 
one can quickly produce generating functions for the recurrences as follows.
Define a sequence of polynomials in
$n, N, x$ by $p_{-1}(N,n,x)=1$, $p_0(N,n,x)=n$, and
\begin{equation}
 p_j(N,n,x):=np_{j-1}(N,n,x) - jn^2(N-j+1)xp_{j-2}(N,n-1,x).
\end{equation}
Then $p_j(N,n,x)$ is a polynomial with coefficients given by the 
coefficients of $a^N_{n-k}$ in
(\ref{eqn_for_recrel}).  The terms of the
recurrence for $a_n^N$ are given by coefficients of $x$ in
$p_N(N,n,x)$, so we have
\begin{equation}
\sum_{k=0}^{\lfloor\frac{N+1}{2}\rfloor}
 c_k^N(n) a^N_{n-k} = 0 ,
\text{ where } p_N(N,n,x)=\sum c_r^N(n) x^r.
\end{equation}

\end{remark}

From (\ref{recrel2}), using the standard method
we obtain the following differential equation
for $f_N(x)=\sum_{n\ge 0} a^N_nx^n,$  where
$\Theta=x\frac{d}{dx}$.
\begin{equation}
\mathcal F_N=
\sum_{k\ge 0}
t^k
\!\!\!\!
\sum_{\scriptsize{\begin{array}{c}
1\le i\le k,
\alpha_i+\beta_i=N+1,\\
\alpha_i\in\mathbb N,\> \alpha_{k+1}=0,\\
1< \alpha_{i+1}+1< \alpha_i\le N,
\end{array}
}}
\!\!\!\!\!\!\!\!\!\!\!\!\!\!\!\!\!\!\!\!\!\!\!\!\!
(\Theta+k)^{N+1-\alpha_1}
\prod_{i=1}^k
-\alpha_i\beta_i\left(\Theta+k-i\right)^{\alpha_i-\alpha_{i+1}}
.
\label{diffeqn1}
\end{equation}
The first few cases (up to sign) are given in Table~\ref{tableofPFensa=1}.

%second table was here

Equations 
$\mathcal F_3$ is the Picard-Fuchs equation of the family of elliptic
curves for $\Gamma_1(6)$, and appears in \cite[Table~7]{SB}.
$\mathcal F_4$ is the Picard-Fuchs
equation of the $A_3$ family of K3 surfaces 
studied in \cite{verrillKyoto}, where $\mathcal F_5$ is also given.
$\mathcal F_5$ and $\mathcal F_6$
 can be found in examples 
\#34 and
\#130 respectively in the tables of
\cite{AZ}.

\begin{table}
\begin{eqnarray*}
\mathcal F_2&=&(4x - 1)\Theta + 2x\\
\mathcal F_3&=&
(9x-1)(x-1)\Theta^2 + 2x(9x - 5)\Theta + 3x(3x - 1)\\
\mathcal F_4&=&(16x-1)(4x-1)
\Theta^3 + 6x(32x - 5)\Theta^2 + 6x(32x^2 - 3)\Theta + 4x(16x - 1)
\\
\mathcal F_5&=&
(25x-1)(9x-1)(x-1)\Theta^4 + 
x(1350x^2 - 1036x + 70)\Theta^3  \\
&&+ x(2925x^2 - 1580x + 63)\Theta^2 + x( 2700x^2 - 1088x + 28)\Theta
\\
&& + 
5x( 180x^2 - 57x + 1)
\\
\mathcal F_6&=&
(4x-1)(16x-1)(36x-1)\Theta^5 + (17280x^3 - 3920x^2 + 140x)\Theta^4\\
&& + (50688x^3 - 8076x^2 + 168x)\Theta^3 + (72576x^3 - 8548x^2 + 112x)\Theta^2\\
&& + (50688x^3 - 4628x^2 + 40x)\Theta + 6x(2304x^2 - 170x + 1)
\\
\mathcal F_7&=&
(49x-1)(25x-1)(9x-1)(x-1)\Theta^6 \\
&&+ (132300x^4 - 116244x^3 + 11844x^2 - 252x)\Theta^5 \\
&&+ (639450x^4 - 431406x^3 + 30798x^2 - 378x)\Theta^4 \\
&& + (1587600x^4 - 844776x^3 + 44232x^2 - 336x)\Theta^3\\
&& + (2127825x^4 - 919770x^3 + 36789x^2 - 180x)\Theta^2\\
&& + (1455300x^4 - 527112x^3 + 16698x^2 - 54x)\Theta \\
&&+ 7x(56700x^3 - 17720x^2 + 459x - 1)
\end{eqnarray*}
\caption{Differential equations for $\sum_{n\ge 0} a_n^Nt^n$ for $2\le N\le 7$,
obtained as special cases of (\ref{diffeqn1})
}
\label{tableofPFensa=1}
\end{table}

\section{Weighted sums}
\label{weightedsumsmethod}

\comment{
Now we consider the equation corresponding to 
$$\left(\sum_{i=1}^N X_i
\right)\left(
 \sum_{i=1}^N\frac{a_i}{X_i}\right) = \lambda.$$
}
Now we want to find a recurrence relation for the terms
\begin{equation}
a_n^{\alpha_1,\dots,\alpha_N}
:=\sum_{\sum p_i=n}
\left(\prod_{i=1}^N
\alpha_i^{p_i}\right)
\left(\frac{n!}{p_1!p_2!\cdots p_N!}\right)^2.
\label{def_gen_an}
\end{equation}
We will call these weighted sums, as opposed to the case where all $\alpha_i=1$.

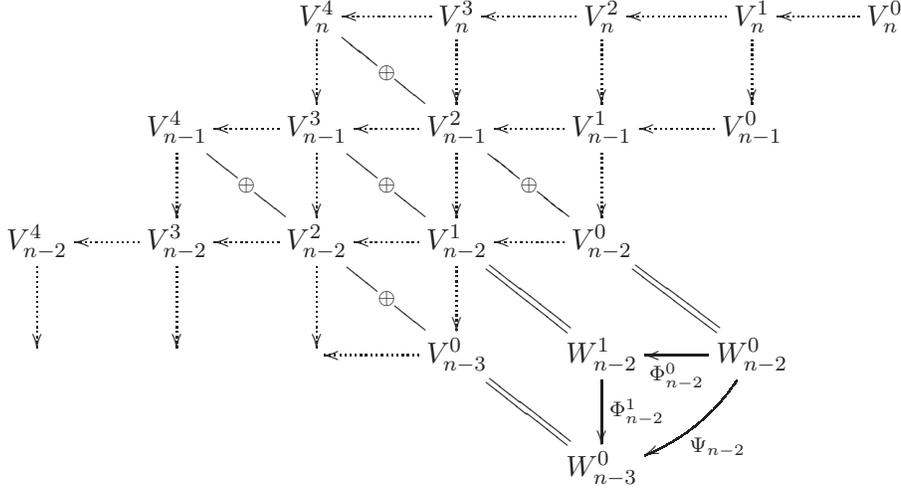
\begin{figure}
$$
\hspace{-1cm}
\xymatrix{
&  &           
 & V_n^4\ard\ar@{-}[dr]|\oplus
 & V_n^3\arl\ard & V_n^2\arl\ard & V_n^1\arl\ard & V_n^0\arl\\ 
  &
& V_{n-1}^4\ard\ar@{-}[dr]|\oplus
 & V_{n-1}^3\arl\ard\ar@{-}[dr]|\oplus
 & V_{n-1}^2\arl\ard
\ar@{-}[dr]|\oplus
 & V_{n-1}^1\arl\ard
 & V_{n-1}^0\arl\\ 
 & V_{n-2}^4\ard & V_{n-2}^3\arl\ard & V_{n-2}^2\arl\ard\ar@{-}[dr]|\oplus
 & 
V_{n-2}^1\arl\ard\ar@{=}[dr]
 & V_{n-2}^0\arl
\ar@{=}[dr]
\\ 
{} & {} & {} & {} & V_{n-3}^0\arl\ar@{=}[dr]
 & W_{n-2}^1\ar[d]^{\Phi_{n-2}^1}
 & W_{n-2}^0\ar[l]^{\Phi_{n-2}^0}\ar@/^3mm/[dl]^{\Psi_{n-2}}
\\
{} & {} & {} & {} & {} & W_{n-3}^0 
}
$$
\caption{Maps used to find the $a_N^{\alpha_1,\dots,\alpha_N}$ recurrence relation}
\label{fig:generaldiagramofmaps}
\end{figure}

Finding a recurrence relation for these $a_n$ is achieved in a similar
manner as for the previous case, but now the situation is slightly more
complicated.  In particular, the version of (\ref{eqn_aux_reln}) 
cannot in this
case be so easily reversed to give a relation like 
(\ref{eqn_aux_reln_reversed})
(this can be done
for $N=2$, but 
for large $N$ the
situation becomes too complicated).
So, although we will obtain a diagram (Figure~\ref{fig:generaldiagramofmaps})
similar to Figure~\ref{fig_arrows1} our
arrows will now be in the opposite direction.  But we will still be able
to obtain a relation using the finite dimensionality of certain 
vector spaces.

Our auxiliary terms in this case are defined by
\begin{equation}
a_n^{\eps_1,\dots,\eps_N}:=
\sum_{\sum p_i=n}
p_1^{\eps_1}\cdots p_N^{\eps_N}
\left(\prod_{i=1}^N
\alpha_i^{p_i}\right)
\left(\frac{n!}{p_1!p_2!\cdots p_N!}\right)^2 
,
\label{eqn_aux_terms_general}
\end{equation}
where $\eps_k\in\{1,-1\}$,
which reduces to (\ref{eqn_aux_terms}) in the case that all $\alpha_i=1$,
with $j$ in (\ref{eqn_aux_terms}) given by $j=\sum \eps_k$.

The generalisation of (\ref{eqn_aux_reln}) is given by
\begin{eqnarray}
na_n^{\eps_1,\dots,\eps_N}
&=&
\sum_{i=1}^N(1-\eps_i)a_n^{\eps_1,\dots,\eps_i+1,\dots,\eps_N}
\nonumber
\\
&&+
n^2\sum_{i=1}^N \alpha_i\eps_i a_{n-1}^{\eps_1,\dots,\eps_i-1,\dots,\eps_N}
\label{eqn_aux_rel_general}
\end{eqnarray}
Notice that since $\eps_i\in\{-1,1\}$,
the right hand side of this expression always has $N$ terms.

The ``direction'' of this equation can not easily be reversed.
Instead of a diagram as in Figure~\ref{fig_arrows1}, built of
figures of the form (\ref{eqn_meaningofdiagram}), we would have
building blocks of a much more complicated form, e.g., such as:

\begin{equation}
\raisebox{4.2ex}{\mbox{
\xymatrix@R=0cm{
a_n^{1,1,0}  & 
{}\save[]+<0cm,-0.3cm>
{\phantom{MmM}\raisebox{-0.6cm}{$a_n^{1,0,0}$}\ar[l]_{1/n}\ar[ld]^{1/n}
\ar[dddd]^{n\alpha_1}}
\restore
\\
a_n^{1,0,1} \\
\\
\\
& a_{n-1}^{0,0,0}
}}}
\end{equation}
meaning $$a_n^{1,0,0}=
 \mbox{$\frac{1}{n}$}a_n^{1,0,1}
+\mbox{$\frac{1}{n}$}a_n^{1,1,0} + n\alpha_1a_{n-1}^{0,0,0}.$$
In order to work with these relations, we consider $N$ fixed, and
define vector spaces
\begin{equation}
V_n^j=\bigoplus_{\sum \eps_i=j}K a_n^{\eps_1,\dots,\eps_N},
\label{eqn_Vndef}
\end{equation}
where $K$ is the field of fractions of $\ZZ[n,\alpha_1,\dots,
\alpha_n]$, and $n$ is considered
an independent variable, and the $a_n$ satisfy only the relations given by
(\ref{eqn_aux_rel_general}).
We now consider (\ref{eqn_aux_rel_general}) as defining a map
$V_n^j\rightarrow V_n^{j+1}\oplus V_{n-1}^{j-1}$.  
Summing such expressions together, we define maps
\begin{equation}
\Phi_n^j:W_n^j\rightarrow W_{n-j}^{1-j},
\text{ where $j=0$ or $1$, and where }
W_n^j:=\bigoplus_{k=0}^{\lfloor\frac{N-1}{2}\rfloor}V_{n-k}^{2k+j}.
\label{eqn_def_Wn}
\end{equation}
We define the composition 
\begin{equation}\Psi_n=\Phi_n^1\circ\Phi_n^0.
\label{def_composistion}
\end{equation}
We can visualise this,
for example in the case of $N=4$ in the diagram in Figure~\ref{fig:generaldiagramofmaps}, 
which shows only part of an infinite diagram.  This diagram corresponds to Figure~\ref{fig_arrows1},
with 
$g(N,4)a_n^{N,5}$ removed, since we will
set $j=N$ so that this is zero, and with
the $a_n^{N,j}$ replaced by $V_n^j$, since in the special case where all $\alpha_i=1$,
the $V_n^j$ would be one dimensional, spanned by $a_n^{N,j}$.

\begin{lemma}
$\dim W_n^j\le 2^{N-1}$.
\label{lem:Wn_dim_bound}
\end{lemma}
\begin{proof}
From the definition (\ref{eqn_Vndef}), $\dim V_n^j\le \binom{N}{j}$
(with equality only if the $a_n^{\bf \eps}$ spanning $V_n$ are
independent over $K$),
and so from the definition (\ref{eqn_def_Wn}), the dimension
 of $W_n^j$ is at most the sum of alternative binomial coefficients.
These sums represent the number of ways of taking an even, respectively 
odd, subset from a set of $N$ items.  Subsets correspond to binary sequences
in $(\ZZ/2\ZZ)^N$, so there is a bijection between these sets corresponding to
the map ${\bf x}\mapsto{\bf x}+(1,0,\dots,0)$.
\end{proof}
\begin{corollary}
For $\alpha_1,\dots,\alpha_N\in\CC$,
if $a_n=a_n^{\alpha_1,\dots,\alpha_N}$ as defined in (\ref{def_gen_an}), then
$a_n,a_{n-1},\dots,a_{n-2^{N-1}}$ satisfy a linear relation over $K$.
\label{corl_lin_rel}
\end{corollary}
\begin{proof}
This is because the maps $\Phi_n^j$ 
(\ref{eqn_def_Wn}) 
are injective, since they
are defined by equalities (\ref{eqn_aux_rel_general}).  A 
linear
relationship (which exists by Lemma~\ref{lem:Wn_dim_bound})
between
$\Psi_{n-2^{N-1}+1}\circ\cdots\circ\Psi_n(a_n),
\Psi_{n-2^{N-1}+1}\circ\cdots\circ\Psi_{n-1}(a_{n-1}),\dots,a_{n-2^{N-1}}$
in $W_{n-2^{N-1}}$ thus pulls back to a relationship between the $a_n$.
\end{proof}
Viewing $a_n$ as numbers rather than independent variables, one obtains
immediately the following result. 

\begin{corollary}
\label{cor:mainresult}
There is a recurrence relation satisfied by the 
$a_n^{\alpha_1,\dots,\alpha_N}$, with at most $2^{N-1}+1$
terms, and with coefficients given as polynomials
in $n$ and $\alpha_1,
\dots,\alpha_N$.
\end{corollary}

\begin{remark}
By making a more careful analysis of the entries of the matrices $\Psi_i$,
it should be possible to obtain a bound  on the degree of the 
coefficients (as polynomials in $n$) of the recurrence relation.
\end{remark}
Note that the recurrence relation may have fewer than $2^{N-1}+1$
terms---by Lemma~\ref{lem:Wn_dim_bound},
this is determined by the dimension of the vector spaces $W_n^j$, which
depends on the relationships between the $\alpha_i$.

\section{Formulae and examples}

Although the above results give a procedure for determining a recurrence relation,
we still do not have a general formula comparable with (\ref{recrel2}).
The recurrence relations must be determined case by case.
For a given set of $\alpha_i$, since we now know a bound on the degree and number
of terms of a recurrence relation, it is a simple matter of linear algebra to determine
the relation explicitly.  
However, what makes this result more interesting is that we can obtain formulas for
an infinite number of cases at once.

\subsection{The recurrence when $N=2$}

In the case $N=2$, 
applying the above method we have the following result.
\begin{theorem}
For $a,b\in\CC$, the terms
$$a_n^{a,b}:=\sum_{p+q=n}a^pb^q\binom{n}{p}^2$$
satisfy a recurrence relation
\begin{equation}
na_n^{a,b} - (2n-1)(a+b)a_{n-1}^{a,b}
 + (n-1)(a-b)^2a_{n-2}^{a,b}=0.
\label{recrel_generalN=2case}
\end{equation}
\label{thm_ab_case}
\end{theorem}
\begin{remark}
Note that substituting $a=b=1$ into (\ref{recrel_generalN=2case})
gives the same result as substituting $N=2$ into (\ref{recrel2}) in
Theorem~\ref{thm:recrel_for_N=1case}.
\end{remark}
\begin{proof}
In this case from (\ref{eqn_aux_rel_general}) we have
\begin{eqnarray*}
na_n^{1,1}&=&n^2aa_{n-1}^{0,1} + n^2ba_{n-1}^{1,0}\\
na_{n}^{0,0}&=&a_{n}^{1,0} + a_{n}^{0,1}\\
na_n^{1,0}&=&a_n^{1,1} + n^2aa_{n-1}^{0,0}\\
na_n^{0,1}&=&a_n^{1,1} + n^2ba_{n-1}^{0,0}.
\end{eqnarray*}
In matrix notation, these are:
\begin{equation}
\left(\begin{matrix}
a_{n}^{0,0}
\\
a_{n+1}^{1,1}
\end{matrix}
\right)
=
\left(\begin{array}{cc}
1/n & 1/n\\
(n+1)b & (n+1)a
\end{array}
\right)
\left(\begin{matrix}
a_{n}^{1,0}
\\ 
a_{n}^{0,1}
\end{matrix}
\right)
\end{equation}
\begin{equation}
\left(\begin{matrix}
a_{n}^{1,0}
\\
a_{n}^{0,1}
\end{matrix}
\right)
=
\left(\begin{array}{cc}
 na&1/n \\
 nb&1/n 
\end{array}
\right)
\left(\begin{matrix}
a_{n-1}^{0,0}
\\ 
a_{n}^{1,1}
\end{matrix}
\right)
\end{equation}
Thus the composition (\ref{def_composistion}) is given by
$$\Psi_n=
\left(\begin{array}{cc}
na & nb\\
1/n & 1/n
\end{array}
\right)
\left(\begin{array}{cc}
1/n & (n+1)b\\
1/n & (n+1)a
\end{array}\right)
=
\left(\begin{array}{cc}
a + b   
& 2abn(n+1)
\\
2/n^2 
&
(a + b)(n+1)/n
\end{array}\right),
$$
and (writing $\ast$ for entries we are not interested in)
$$
\Psi_{n-1}=
\left(\begin{array}{cc}
a + b   
& \ast
\\
2/(n-1)^2 
&
\ast
\end{array}\right),
\>\>
\Psi_{n-1}\circ\Psi_n=
\left(\begin{array}{cc}
(a+b)^2 + 4ab(n-1)/n
& \ast %
\\
\frac{2(2n-1)(a+b)}{(n-1)^2n}
&
\ast
\end{array}\right).
$$
Now we can see that the linear relationship referred to in
Corollary~\ref{corl_lin_rel}, which now is a relationship in
$W_{n-2}$ spanned by $a^{0,0}_{n-2}$ and $a^{1,1}_{n-1}$,
is given by
$$
n
\left(\begin{matrix}
(a+b)^2 + \frac{4ab(n-1)}{n}
\\
\frac{2(2n-1)(a+b)}{(n-1)^2n}
\end{matrix}\right)
-(2n-1)(a+b)
\left(\begin{matrix}
a + b   
\\
\frac{2}{(n-1)^2}
\end{matrix}\right)
+
(n-1)(a-b)^2
\left(\begin{matrix}
1
\\
0
\end{matrix}\right)=0.
$$
The recurrence relation follows, as
described in Corollary~\ref{corl_lin_rel}.
\end{proof}
\begin{corollary}
The function $f=\sum_{n\ge 0}a_n^{a,b}t^n$ satisfies the following
differential equation:
$$\left[(1 - 2(a+b)t+t^2(a-b)^2)\Theta - t(a+b)+t^2(a-b)^2\right]
f=0.$$
\end{corollary}
\begin{remark}
Notice that the coefficient of $\Theta$ factors as 
$((\sqrt{a} + \sqrt{b})^2t-1)((\sqrt{a} - \sqrt{b})^2t-1).$
This is an example of a general phenomena, which has a geometric explanation,
given by the determination of the singular members of the family 
$\mathcal X^{\alpha_1,\dots,\alpha_N}_t$; the $A_4$ situation is given in \cite[Lemma~3.7]{HV},
and the general case works in exactly the same way.
\end{remark}

\begin{example}
\label{ex:11_1-1}
For $a=b=1$ in (\ref{recrel_generalN=2case}) in Theorem~\ref{thm_ab_case}
we obtain the recurrence
\begin{equation}
na_n^{1,1} = 2(2n-1)a_{n-1}^{1,1}.
\end{equation}
In the introduction we remarked that $a_n^{1,1}=\binom{2n}{n}$.
In the case $a=-b=1$ we have
\begin{equation}
\label{eq:recrel1-1}
na_n^{1,-1} = 4(n-1)a_{n-2}^{1,-1}.
\end{equation}
From the definition (\ref{eqn:firstdefinition})
 we see that $a_1^{1,-1}
=0$, so from (\ref{eq:recrel1-1}) all odd terms are zero.
Setting $c_n=a^{1,-1}_{2n}$ and substituting in (\ref{eq:recrel1-1})
we obtain $2na_2n^{1,-1} = 4(2n-1)a_{2n-2}^{1,-1}$, and then
\begin{equation}
\label{eq:recrel1-1c}
nc_n = 2(2n-1)c_n
\end{equation}
Since this is the same relation as (\ref{ex:11_1-1}), and
$a_0^{1,1}=c_0=1$, we have $a_n=c_n$, giving us the formula 
(\ref{motivating_1}) given in the introduction.
Note, both these examples are well known.  Refer also to
\cite{C} and \cite{}.
\end{example}

\subsection{The recurrence when $N=3$}

In the case of $N=3$, the matrices we are interested in are given by
$$
\Phi_n^0=
\left(
\begin{array}{cccc}
0       & 1/n     & 1/n  & 1/n       \\
1/(n+1) & 0       & (n+1)c & (n+1)b\\
1/(n+1) & c(n+1) & 0       & (n+1)a\\
1/(n+1) & b(n+1) & (n+1)a & 0
\end{array}
\right)^t
$$
and
$$
\Phi_n^1=
\left(
\begin{array}{cccc}
0   & (n+1)a&  (n+1)b&  (n+1)c\\
na & 0      &  1/n    & 1/n\\
nb & 1/n    &  0      & 1/n\\
nc & 1/n    &  1/n    & 0
\end{array}
\right)^t.
$$
From this, setting $A=a+b+c$, we obtain
$$
\Psi_n=
\left(
\begin{array}{cccc}
A & 2cbn(n+1)   & 2acn(n+1) & 2ban(n+1)
\\
2/n^2       & A+(b + c)/n & 2a + a/n & 2a + a/n
\\
2/n^2       & 2b + b/n & A + (a + c)/n & 2b + b/n
\\
2/n^2       & 2c + c/n & 2c + c/n & A + (a + b)/n
\end{array}
\right)^t.
$$
In $W_{n-4}$,
$\Psi_{n-3}\circ\cdots\circ\Psi_n(a_n),\dots,a_{n-4}$
are given by the columns of the following matrix (computer algebra packages e.g.,
\cite{pari}, \cite{magma},
are helpful in obtaining and manipulating this matrix).
\begin{equation}
\label{bigmatirxN=3}
\left[
\begin{matrix}
H_0
&H_1 &A^2+\frac{4B(n-3)}{n-2} &A&1\\
G_2(a,b,c) &G_1(a,b,c)&G_0(a,b,c,n-2) &\frac{2}{(n-3)^2}&0\\
G_2(b,c,a)&G_1(b,c,a)&G_0(b,c,a,n-2) &\frac{2}{(n-3)^2}&0\\
G_2(c,a,b)&G_1(c,a,b)&G_0(c,a,b,n-2) &\frac{2}{(n-3)^2}&0
\end{matrix}
\right]
\end{equation}
where 
$A=a + b + c$, $B=ab+bc+ca$, and
\begin{eqnarray*}
H_0&=&
A^4 + \frac{4A^2B(6n^3 - 28n^2 + 37n - 12)}{n(n-1)(n-2)}
+\frac{16B^2(n-1)(n-3)}{n(n-2)}
\\
&&+ \frac{4Aabc(10n^2-10n+3)(4n-7)(n-3)}{n^2(n-1)^2}
+\frac{12Aabc(4n-3)(2n-5)}{n^2},
\\
H_1&=&A^3 + \frac{4(3n^2 - 13n + 13)AB}{(n-1)(n-2)}
+\frac{12(n - 3)(4n - 7)abc}{(n-1)^2},
\end{eqnarray*}
$$
\begin{array}{l}
\!\!\!\!\!\!
G_0(a,b,c,n)=
\frac{2(n-1)(4n-3)a+2An(2n-1)}{n^2(n-1)^2},\\
G_1(a,b,c)=
\frac{(4n-7)(10A-4a)a}{(n-1)^2(n-3)}
-\frac{6(2n-3)Aa               }{(n-1)^2(n-2)^2(n-3)}
-\frac{2(n-2)(A^2-4B)               }{(n-1)(n-3)^2}
+\frac{2(2n-5)(2n-3)A^2}{(n-1)(n-2)(n-3)^2},
\end{array}
$$
$$
\begin{array}{l}
G_2(a,b,c)=
\Big[A^2+\frac{4B(n-1)}{n}\Big]\frac{(2(n-3)(4n-11)a+2A(n-2)(2n-5))}{(n-2)^2(n-3)^2}
\\
-A(2n-1)\frac{(n-3)4(6n^2-22n+19)bc-2(n-3)(12n^2-50n+47)a^2 - 2B(2n-5)^2(5n-7)}
{n(n-1)^2(n-3)^2(n-2)}
\\
+a(4n-3)\frac{4cb(12n^3-71n^2+138n-87)-2(n-3)(2n-3)a^2+8(a^2+b^2+c^2)(n-2)^2(n-3)+2B(n - 3)(2n - 5)(4n - 5)}{(n-3)^2(n-2)(n-1)n^2}
\\
-2\Big[(4n-3)(n-1)a+An(2n-1)\Big]\frac{(2n-3)(2n-5)(a^2 + cb) + (b+c)a(4n-7) + A^2(-n^2 + 3n -2)}
{(n - 3)(n - 2)(n - 1)^2n^2}.
\end{array}
$$
In general (\ref{bigmatirxN=3}) has rank $4$, which can be verified by
plugging in any values for $a,b,c$.  If exactly two of $a,b,c$ are equal, then
two of the last three rows are equal, so the rank is $3$; if $a=b=c$ the rank is $2$.
By the same method as for Theorem~\ref{thm_ab_case}
a recurrence relation is obtained by finding linear relations between the columns, i.e.,
by finding the kernel of this matrix.  Generally this gives a $5$ term recurrence,
but if $a=b$, or $a=b=c$, there will be a relation between the last $4$ or $3$ columns
respectively.
By explicitly computing the kernel (again, aided by \cite{magma} and \cite{pari})
 we obtain
the following result.
\begin{theorem}
\label{thm_abc_case}
For $a,b,c\in\CC$, $abc\not=0$ the terms
\begin{equation}
a_n=a_n^{a,b,c}:=\sum_{p+q+r=n}a^pb^qc^r\binom{n}{p,q,r}^2
\label{def:A3caseterms}
\end{equation}
satisfy a recurrence relation $R(n)=0$ where
\begin{eqnarray}
\nonumber
R(n):=F_{11}F_7 n^2 a_n 
-AF_{11}\Big[
2(n-1)(2n-3)(4n-1) + 3\Big]a_{n-1}\\
\nonumber
+\Big[(2n - 3)^2F_{11}F_3A^2  + \left[F_3F_{11}(n - 2)(n - 1) - 3\right]s
\Big]a_{n-2}\\
\label{recrelN=3abc}
-F_9F_3\Big[sA(4n^2-18n+19) + 4abcF_{11}F_7\Big]a_{n-3}
+F_7F_3s^2(n-3)^2a_{n-4},\phantom{MM}
\end{eqnarray}
and $F_{p}=4n-p, A=a+b+c$ and
$s=a^2 + b^2 + c^2 - 2ab - 2ac - 2cb$.
\end{theorem}
\begin{remark}
The differential equation for $\sum a_nt^n$
 corresponding to (\ref{recrelN=3abc}) has degree $4$ in 
$\Theta=\frac{td}{dt}$; the coefficient of $\Theta^4$ is
$\prod((\sqrt{a}\pm\sqrt{b}\pm\sqrt{c})^2t-1)$.
\end{remark}

\begin{example}
If $s=a^2 + b^2 + c^2 - 2ab - 2ac - 2cb=0$, which happens when $\sqrt{a}\pm\sqrt{b}\pm\sqrt{c}=0$,
 then (\ref{recrelN=3abc}) is
divisible by $F_{11}$. Since $n$ is an integer, so
$F_{11}$ is never $0$, we obtain
\begin{eqnarray*}
0=F_7 n^2 a_n -A\left[2(n-1)(2n-3)(4n-1) + 3\right]a_{n-1}\\
+A^2(2n - 3)^2F_3a_{n-2}
-4abcF_9F_3F_7a_{n-3}.
\end{eqnarray*}
\end{example}
\begin{example}
\label{ex:omega3}
If we have $A=s=0$, which happens (up to permutations
and scaling) for $(a:b:c)=(1:\omega:\omega^2)$, 
where $\omega$ is a primitive root of $1$, we
obtain
\begin{equation}
a_n:=\sum_{p+q+r=n}\omega^{p-q}\binom{n}{p,q,r}^2
\>\>\>\Rightarrow\>\>\>
n^2 a_n =4 (4n-9)(4n-3)a_{n-3}.
\end{equation}
The corresponding
differential equation for $\sum a_nt^n$ is  
\begin{equation}
(1-64t^3)\Theta^2 + 64t^3\Theta + 20.
\end{equation}
Note that though the $a_n$ are defined using $\omega$, all terms of the sequence are
integers (since $\sum_{i=1}^3 \omega^{ni}=0$ or $1$ depending on $n\mod 3$).
The first few terms are:
$$    1,
    0,
    0,
    12,
    0,
    0,
    420,
    0,
    0,
    18480,
    0,
    0,
    900900,\dots$$
(This is sequence A000897 in Sloane's table of integer sequences 
\cite{sloan}.)
Since the recurrence relation is so simple it is easy to verify 
that the nonzero terms are
given by  $\binom{4m}{2m,m,m}$, where $m=n/3$.
(For example, show that $\binom{4m}{2m,m,m}$ satisfies the above
relation; given $a_0=1$ the solution is unique.)
In other words,
\begin{equation}
\sum_{p+q+r=3m}\omega^{p-q}\binom{3m}{p,q,r}^2
=\binom{4m}{2m,m,m}.
\end{equation}

\end{example}

\subsection{Examples when $N=4$}
Finding explicit equations for the general case for  higher $N$ is possible but 
time consuming,
so for $N=4$ 
we just give a couple of examples, in cases where the recurrence has fewer terms
than in general (when there will be $9$ terms), due to relationships between the $\alpha_i$.

\begin{example}
For $\{\alpha_i\}=(1,1,1,9)$ the first few $a_n$ are
$$ 1, 12, 204, 4224, 99324, 2546352, 69359424, 1973611008, 58005903708,\dots,$$
and the
recurrence relation obtained by the above method is
\begin{eqnarray*}
(n - 1)n^3(10n^2 - 35n + 31)a_{n}&&\\
-4(n - 1)(140n^5 - 700n^4 + 1289n^3 - 1104n^2 + 477n - 84)a_{n-1}&&\\
+4(1960n^6 - 14700n^5 + 44986n^4 - 71829n^3 + 63127n^2 - 29022n + 5496)a_{n-2}&&\\
-1152(n - 2)^3(2n - 3)(10n^2 - 15n + 6)a_{n-3}\\=0
\end{eqnarray*}
The corresponding differential equation is:
\begin{eqnarray*}
10 (4t - 1)(16t - 1)(36t - 1)\Theta^6
+45 (4608t^3 - 784t^2 + 1)\Theta^5\\
+2  (378432t^3 - 31172t^2 - 222t - 33)\Theta^4
+  (1433088t^3 - 54636t^2 - 148t + 31)\Theta^3\\
+4  t(371520t^2 - 6217t + 36)\Theta^2
+24 t(33408t^2 - 235t + 3)\Theta
+576t^2(306t - 1).
\end{eqnarray*}
\end{example}

\begin{example}
For $\{\alpha_i\}=(1,1,1,-3)$ the first few $a_n$ are
$$ 1, 0, -12, -96, -180, 5760, 70080, 161280, -5144580, -68974080,\dots,$$
and we obtain
\begin{eqnarray*}
(n - 1)n^3(14n^3 - 84n^2 + 165n - 107)a_{n}&&\\
-4(n - 1)^3(2n - 5)(14n^3 - 42n^2 + 39n - 12)a_{n-1}&&\\
+4(14n^3 - 84n^2 + 165n - 107)(28n^4 - 112n^3 + 163n^2 - 102n + 24)a_{n-2}&&\\
+192(n - 2)^2(28n^5 - 210n^4 + 582n^3 - 737n^2 + 426n - 93)a_{n-3}&&\\
+2304(n - 3)^2(n - 2)^2(14n^3 - 42n^2 + 39n - 12)a_{n-4}&&\\=0.
\end{eqnarray*}
\end{example}

\begin{example}
For $\{\alpha_i\}=(1,1,-1,-1)$ the first few $a_n$ are
$$ 1, 0, -4, 0, 156, 0, -5440, 0, 239260, 0, -11151504, 0, 551724096, 0, \dots,$$
and we have
\begin{eqnarray*}
(n - 1)n^3(10n^2 - 55n + 76)a_{n}\\
+4(120n^6 - 1140n^5 + 4282n^4 - 8107n^3 + 8170n^2 - 4176n + 864)a_{n-2}\\
-1024(n - 3)^3(n - 2)(10n^2 - 15n + 6)a_{n-4}=0
\end{eqnarray*}
\end{example}

\begin{example}
The $N=4$ version of Example~\ref{ex:omega3} is
$\{\alpha_i\}=\{1,i,-1,-i\}$.  The first few $a_n$ are
$$
1, 0, 0, 0, -132, 0, 0, 0, 113820, 0, 0, 0, -140078400,0, 0, 0, 201740158620, \dots,
$$
and there is a three term recurrence
\begin{eqnarray*}
0=(n - 3)(n - 2)(n - 1)n^3(48n^4 - 1032n^3 + 8276n^2 - 29347n + 38840)a_{n}\\
+
16(6528n^{10} - 218688n^9 + 3180512n^8 - 26345016n^7 + 137020240n^6 \phantom{MM}\\
 - 465036692n^5 + 1036364052n^4 - 1486439881n^3\phantom{MM}\\
 + 1303139340n^2 - 627480000n + 127008000)
a_{n-4}\\
\hspace{-0.3cm}
+ 2^{12}(n - 7)^2(n - 6)(n - 5)^2(n - 4)(48n^4 - 264n^3 + 500n^2 - 387n + 108)a_{n-8}.
\end{eqnarray*}
The fact that in this expression the two
factors of the form $48n^4+\cdots$
are related by the change of variables $n\mapsto n-4$ leads one to
expect that there is probably a recurrence with more terms but coefficients
of lower degree.
\end{example}

\section{Picard-Fuchs equations in the elliptic curve case}
\label{sec:ellipticPF}

If a sequence satisfies a recurrence relation, then it satisfies
infinitely many.  So, we do not claim that the
recurrence relation obtained by this method give Picard-Fuchs equation.  
However they could be determined from these relations,
once we know (from the geometry) what order the equation should have.
  For example, for the $A_n$ case, the family of elliptic
curves in $\PP^2$ is
\begin{equation}
E_t: (X+Y+Z)(aXY+bYZ+cZX)t=XYZ.
\label{eqn:Eabc}
\end{equation}
Aided by computer algebra \cite{magma}, we find that the $j$-invariant of $E_{1/t}$ is
\begin{eqnarray}
j(E_{1/t})=\frac{
\left(
P
+16abct\right)^3
}
{(abct)^2P
},
\end{eqnarray}
and that
for a certain Weierstrass form $y^2 = x^3 - 27c_4x - 54c_6$
\begin{eqnarray}
c_4(E_{1/t})&=&
\frac{16(a - t)^4}{a^4(b-c)^4}
(P+16abct)\\
c_6(E_{1/t})&=&
\frac{-64(a - t)^6}{a^6(b-c)^6}
\sqrt{\left(P+64abct\right)}
(P-8abct)\\
\text{ where }
P&:=&
\displaystyle{\prod_{\eps_1,\eps_2\in\{1,-1\}}}
\Big(\left(\sqrt{a}+\eps_1\sqrt{b}+\eps_2\sqrt{c}\right)^2 - t\Big).
\end{eqnarray}
A period of this family given (up to a constant factor)
by the series (\ref{def:A3caseterms}), and satisfies
a differential equation corresponding to the recurrence relation (\ref{recrelN=3abc}), though this
is not the Picard-Fuchs equation, since the degree is too high.
However, one can easily show (again a computer algebra package is helpful) that
with $R(n)$ as in (\ref{recrelN=3abc}), the expression
\begin{eqnarray}
\nonumber
&&3(4n-15)R(n) - 2(a+b+c)(4n-11)R(n-1) - 
s(4n-7)R(n-2),\\
&&\text{where }\> s=a^2+b^2+c^2-2(ab+bc+ca),
\end{eqnarray}
is divisible by $(4n-7)(4n-11)(4n-15)$.  After dividing we have a recurrence relation
\begin{eqnarray}
\nonumber
0&=&3n^2a_n - A(14n^2 -18n+7)a_{n-1}\\
\nonumber
&&+\left[A^2(20n^2 - 56n + 41)+s(5n^2 - 6n - 2)\right]a_{n-2}\\
\nonumber
&&-\left[2A^3(2n-5)^2 + 4Asn(3n-7) + 12abcF_3F_9\right]a_{n-3}\\
\nonumber
&&+\left[A^2s(4n^2-16n+9) + 8abcAF_{11}F_{13} + s^2(n^2 + 4n - 23)\right]a_{n-4}\\
&&+\left[As^2(2n^2 - 22n + 57) + abcsF_{17}F_{19}\right]a_{n-5},
\end{eqnarray}
where $F_{p}=4n-p, A=a+b+c$ and
$s=a^2 + b^2 + c^2 - 2ab - 2ac - 2cb$.
The corresponding Picard-Fuchs differential equation 
(assuming $a,b,c$ are distinct) is then
\begin{eqnarray}
\Theta^2 + t\sum_{i=1}^6\frac{\eps_i}{t-u_i}\Theta + 
t^2\sum_{i=0}^6\frac{\beta_i}{t-u_i},
\label{ellipticPFeqn}
\end{eqnarray}
where $u_0=0$, $u_1,\dots,u_4$ are the values of $1/(\sqrt{a}\pm\sqrt{b}\pm\sqrt{c})^2$, 
which are real singularities, corresponding to singular elliptic curves,
and
$u_5$ and $u_6$ 
are apparent singularities, given by the roots of
$st^2 +2At-3=0.$
The denominators are $\eps_1=\dots=\eps_4=1$, $\eps_5=\eps_6=-1$,
$\beta_0=-a-b-c$, $\beta_5=-\frac{3}{4u_5}$,
$\beta_6=-\frac{3}{4u_6}$,
and
\begin{equation}
\beta_i=\frac{1}{32u_i}\left(\eps^1_i\eps^2_i\sqrt{\frac{1}{abcu_i}}{\left(A-\frac{1}{u_i}\right)} + 26\right) 
\text{ for } 1\le i\le 4,
\end{equation}
where $u_i=1/(\sqrt{a}+\eps^1_i\sqrt{b}+\eps^2_i\sqrt{c})^2$.

An alternative method that could be used to find  the Picard-Fuchs equation is
described by Stienstra and Beukers \cite[\S11]{SB}.  
The values of $\eps_i$ and $\beta_i$ could be determined as in
\cite{SB}, by using knowledge of the solutions, which have the form given in \cite[\S11]{SB}, and
the relations given in \cite[\S15.4,16.4]{ince}.

\begin{example}
The Picard-Fuchs equation 
of the family of elliptic curves given by (\ref{eqn:Eabc})
with $a,b,c=1,16,64$ is
\begin{eqnarray*}
0&=&f'' 
+\left[
\begin{array}{r}
\frac{1}{t}
+\frac{9}{9t-1}
+\frac{25}{25t-1}
+\frac{121}{121t-1}
+\frac{169}{169t-1}
-\frac{11} {11t+1} 
-\frac{165}{165t-1}
\end{array}
\right]f'\\
&&+
3\left[
\begin{array}{l}
\frac{-27}{t}
+\frac{27\cdot 131}{2^7(9t-1)}
+\frac{5^4\cdot23}{2^7(25t-1)}
+\frac{11^4\cdot 53}{2^7(121t-1)}\\
-\frac{13^5}{2^7(169t-1)}
+\frac{11^2} {4(11t+1)} 
-\frac{65^2}{4(65t-1)}.
\end{array}
\right]f.
\end{eqnarray*}

\end{example}

\section{Further recurrences}
\label{sec:furtherapplications}
We expect that the method of finding recurrences described in this
paper should be applicable in many other examples.  A simple generalization
would be to replace squares by higher powers.  The simplest case, when
$N=2$, is the sequence $\{c_n^{a,b,k}\}_n$ with
\begin{equation}
c^{a,b,k}_n := \sum_{p=0}^ka^pb^{n-p}\binom{n}{p}^k.
\label{eqn:k=3case}
\end{equation}
In this case we would introduce auxillary terms
$$c_n^{a,b,k,i,j}:=
\sum_{p+q=n}^ka^pb^qp^iq^j\binom{n}{p}^k,\>\>\>k>i, j\ge 0.$$
%For simplicity of notation we will drop 
%$k$ from the superscripts.
In the case that $a=b$ we may take $i\ge j$.
Similarly to Section~\ref{weightedsumsmethod}, we define vector spaces
$V_{k,n}^m$ spanned by $c_n^{a,b,k,i,j}$, where $i+j=m$, for 
$0\le m\le 2(k-1)$,
and we set 
$$W_{k,n}^j:=V_{k,n}^j\oplus V_{k,n+1}^{k+j}\>
0\le j<k;\>\>\>\>
W_{k,n}^k:=V_{k,n}^k.
$$
Now instead of two maps $\Phi_n^0$ and $\Phi_n^1$
(as in (\ref{eqn_def_Wn})), we would have
$k$ maps $\Phi_n^i$ for $i=0,\dots,k-1$, and we would consider a certain
composition $\Psi_n$ of these
(with image in $W_{k,n-1}^0$ (or $W_{k,n-1}^{k-1}$ in the special case
indicated below)),
in order to obtain a relation in one of the $W$ spaces with minimal
dimension.  (It is helpful to draw a diagram of these spaces,
similar to Figure~\ref{fig:generaldiagramofmaps}.)

In the case where $a=b=1$ this seems to be more or less the same method
as used by \cite{C}, and so this method may be viewed as a generalization
of his method.  Generally $\dim V_{k,n}^m=\min(m+1,2k-m-1)$,
and $\dim W_{k,n}^j=k$,
so we would obtain recurrences with $k+1$ terms.
But in the case $a=b=1$, we may identify $c_n^{a,b,k,i,j}$
and $c_n^{a,b,k,j,i}$, so 
%$\dim V_{k,n}^m=\lceil\min(m+1,2k-s-1)/2\rceil$, and 
$\dim W_{k,n}^m=\lfloor\frac{k}{2}\rfloor+1$, unless $k$ is even and $j=k-1$,
in which case $\dim W_{k,n}^{k-1}=\frac{k}{2}$, so we obtain recurrances with
$\lfloor\frac{k+3}{2}\rfloor$ terms, as in  \cite{C}.

\begin{example}
In the case $k=3$ in (\ref{eqn:k=3case}) we have
maps $\Phi_n^0:W_n^0\rightarrow W_n^1$,
$\Phi_n^1:W_n^1\rightarrow W_n^2$,
$\Phi_n^2:W_n^2\rightarrow W_{n-1}^0$
given by the following three
matrices respectively,
$$
\left(
\begin{matrix}
\frac{1}{n} & b(1+n)^2  & 0\\
\frac{1}{n} & 0 &         a(n+1)^2\\
0 & \frac{1}{(n+1)}& \frac{1}{(n+1)}
\end{matrix}
\right)
,
\>\>\>
\left(
\begin{matrix}
\frac{1}{n} & 0 & b(1+n)^2\\
0 & \frac{1}{n} & a(1+n)^2\\
\frac{1}{n}  &\frac{1}{n}  &0
\end{matrix}
\right)
,
\>\>\>
\left(
\begin{matrix}
n^2a & n^2b &0 \\
0 & \frac{1}{n} & \frac{1}{n}\\
\frac{1}{n} &0 &\frac{1}{n}
\end{matrix}
\right)
$$
with respect to the bases 
$\{a_n^{0,0}, a_{n+1}^{1,2}, a_{n+1}^{2,1}\}$ for
$W_n^0$; 
$\{a_n^{1,0}, a_n^{0,1}, a_{n+1}^{2,2}\}$ for
$W_n^1$;
$\{a_n^{2,1}, a_n^{0,2}, a_n^{1,1}\}$ for
$W_n^2$,
and
$\{a_{n-1}^{0,0}, a_{n}^{1,2}, a_{n}^{2,1}\}$ for
$W_{n-1}^0$.
We set $\Psi_n=\Phi_n^2\Phi_n^1\Phi_n^0$, and consider 
$\Psi_{n-2}\Psi_{n-1}\Psi_{n}(a_n),
\Psi_{n-2}\Psi_{n-1}(a_{n-1}),
\Psi_{n-2}(a_{n-2})$ and $a_{n-3}$, given in $W_{n-2}^0$
(with basis $\{a_{n-2}^{0,0}, a_{n-1}^{1,2}, a_{n-1}^{2,1}\}$) by the
columns of the following matrix:
$$
\left(
\begin{matrix}
\frac{3abA(3(n-2)(n-1)^2(9n-4)+4n)}{(n-1)^2n^2}
+A^3&
\frac{6(n - 2)(3n-5)ab}{(n-1)^2} + A^2 & A &1\\
\frac{\ast}{(n-2)^3(n-1)^2n^2}
&\frac{3(n - 1)(2n-3)b + 3a(5n^2 - 16n + 13)}{(n - 2)^3(n - 1)^2} &
\frac{3}{(n-1)^3} &0\\
\frac{\ast}{(n-2)^3(n-1)^2n^2}
&\frac{3(n - 1)(2n-3)a + 3b(5n^2 - 16n + 13)}{(n - 2)^3(n - 1)^2} & 
\frac{3}{(n-1)^3} &0\\
\end{matrix}
\right),
$$
where $A=a+b$, and $\ast$ denoted unenlightening polynomials of
degree $4$ in $n$ and degreee $2$ in $a$ and $b$.
Note that the determinant of the matrix consisting of the
last three columns is
$\frac{9(a+b)^2(a-b)(3n-2)}{n^2(n-1)^2(n-2)^3}$, and so in the case
$a=\pm b$ there are recurrences with at most $3$ terms
(these can be found for example in \cite{W}, so we do not give them here).
Otherwise we expect recurrences with $4$ terms.
In general by finding the kernel of this matrix we find that the sequence
$c_n:=c_n^{a,b,3}$ satisfies a relation
\begin{eqnarray}
\nonumber
0&=&3n^2(3n-5)c_n - (27n^3 - 72n^2 + 51n - 12)(a+b)c_{n-1}\\
\nonumber
&-&\left[
(a+b)^2+(3n-5)\Big(
(a+b)^2+(9ab-(a+b)^2)(3n-4)(3n-2)\Big)\right]c_{n-2}\\
&-&3(a+b)^3(n-2)^2(3n-2)c_{n-3}.
\end{eqnarray}

\end{example}

\noindent
Helena A. Verrill,\\
Department of Mathematics,\\
Louisiana State University\\
Baton Rouge, LA 70803-4918, USA\\
{\tt verrill@math.lsu.edu}

\end{document}